\documentclass[12pt]{article}
\usepackage{amsmath,amssymb,amsthm,comment}

\newtheorem{theorem}{Theorem}[section]
\newtheorem{thmy}{Theorem}

\def\barr{\begin{array}}
\def\earr{\end{array}}

\title{Finite groups with dense ${\cal CD}$-subgroups}
\author{Ryan McCulloch and Marius T\u arn\u auceanu}
\date{December 30, 2024}

\begin{document}

\maketitle

\begin{abstract}
A group $G$ is said to have dense ${\cal CD}$-subgroups if each non-empty open interval of the subgroup lattice $L(G)$ contains a subgroup in the Chermak--Delgado lattice ${\cal CD}(G)$. In this note, we study finite groups satisfying this property.
\end{abstract}

{\small
\noindent
{\bf MSC2000\,:} Primary 20D30; Secondary 20D60, 20D99.

\noindent
{\bf Key words\,:} finite group, Chermak--Delgado lattice, subgroup lattice, density.}

\section{Introduction}

Let $\mathcal{X}$ be a property pertaining to subgroups of a group. We say that a group $G$ has \textit{dense $\mathcal{X}$-subgroups} if for each pair $(H,K)$ of subgroups of $G$ such that $H<K$ and $H$ is not maximal in $K$, there exists a $\mathcal{X}$-subgroup $X$ of $G$ such that $H<X<K$. Groups with dense $\mathcal{X}$-subgroups have been studied for many properties $\mathcal{X}$: being a normal subgroup \cite{9}, a pronormal subgroup \cite{19}, a normal-by-finite subgroup \cite{5}, a nearly normal subgroup \cite{6} or a subnormal subgroup \cite{7}.

In our note, we will focus on the property of being a ${\cal CD}$-subgroup, that is a subgroup contained in the Chermak--Delgado lattice. Let $G$ be a finite group and $L(G)$ be the subgroup lattice of $G$. We recall that the \textit{Chermak--Delgado measure} of a subgroup $H$ of $G$ is defined by
\begin{equation}
m_G(H)=|H||C_G(H)|.\nonumber
\end{equation}Let
\begin{equation}
m^*(G)={\rm max}\{m_G(H)\mid H\leq G\} \mbox{ and } {\cal CD}(G)=\{H\leq G\mid m_G(H)=m^*(G)\}.\nonumber
\end{equation} Then the set ${\cal CD}(G)$ forms a modular, self-dual sublattice of $L(G)$, which is called the \textit{Chermak--Delgado lattice} of $G$. It was first introduced by Chermak and Delgado \cite{4}, and revisited by Isaacs \cite{8}. In the last years there has been a growing interest in understanding this lattice (see e.g. \cite{1}-\cite{3}, \cite{10}-\cite{14}, \cite{16}-\cite{18} and \cite{20}). We also recall several important properties of the Chermak--Delgado measure:
\begin{itemize}
\item[$\bullet$] if $H\leq G$ then $m_G(H)\leq m_G(C_G(H))$, and if the measures are equal, then $C_G(C_G(H))=H$;
\item[$\bullet$] if $H\in {\cal CD}(G)$ then $C_G(H)\in {\cal CD}(G)$ and $C_G(C_G(H))=H$;
\item[$\bullet$] the maximal member $M$ of ${\cal CD}(G)$ is characteristic and satisfies ${\cal CD}(M)={\cal CD}(G)$;
\item[$\bullet$] the minimal member $M(G)$ of ${\cal CD}(G)$ (called the \textit{Chermak--Delgado subgroup} of $G$) is characteristic, abelian and contains $Z(G)$.
\end{itemize}

Our main results are stated as follows. The first one gives some information about finite $p$-groups having dense ${\cal CD}$-subgroups. Obviously, finite $p$-groups of order $p$ satisfy this property.

\begin{theorem}
Let $G$ be a finite $p$-group of order $p^n$, $n\geq 2$, having dense ${\cal CD}$-subgroups. Then:
\begin{itemize}
\item[{\rm a)}] $|Z(G)|=p$ and $m^*(G)=p^{n+1}$.
\item[{\rm b)}] All subgroups of order $p^2$ of $G$ contain $Z(G)$. Moreover, all normal subgroups of order $p^2$ of $G$ are contained in ${\cal CD}(G)$.
\item[{\rm c)}] ${\rm Im}(m_G)=\{p^n,p^{n+1}\}$.
\end{itemize}
\end{theorem}

Theorem 1.1 shows that the finite $p$-groups with dense ${\cal CD}$-subgroups are a subclass of the class $\cal C$ described in Section 3 of \cite{17}. Note that all extraspecial $p$-groups belong to class $\cal C$, but only some of the extraspecial $p$-groups have dense ${\cal CD}$-subgroups.  For $p$ odd, all extraspecial groups of order $p^3$ have dense ${\cal CD}$-subgroups, and all extraspecial groups of order $p^k$, $k \geq 5$, do not have dense ${\cal CD}$-subgroups (such groups have subgroups of order $p^2$ that intersect the center trivially).  For $p=2$ the situation is similar, except that one of the extraspecial groups of order $32$ has dense ${\cal CD}$-subgroups and the other one does not.
\bigskip

The second result classifies finite groups which are not $p$-groups, but have dense ${\cal CD}$-subgroups.

\begin{theorem}
Let $G$ be a finite group with $|\pi(G)|\geq 2$.  The group $G$ has dense ${\cal CD}$-subgroups if and only if $G$ is a nonabelian group of order $pq$, where $p$ and $q$ are primes.
\end{theorem}

Most of our notation is standard and will not be repeated here. Basic definitions and results on groups can be found in \cite{8}. For subgroup lattice concepts we refer the reader to \cite{15}.

\section{Proofs of the main results}

\noindent{\bf Proof of Theorem 1.1.} 
\begin{itemize}
\item[{\rm a)}] Let $H\leq G$ with $|H|=p^2$. By hypothesis, there is $K\in{\cal CD}(G)$ such that $1<K<H$. Then $|K|=p$ and since $Z(G)\subseteq K$, it follows that $Z(G)=K$. Thus $|Z(G)|=p$, $Z(G)\in{\cal CD}(G)$ and $m^*(G)=p^{n+1}$.
\item[{\rm b)}] Since every subgroup of order $p^2$ of $G$ contains a ${\cal CD}$-subgroup of order $p$, we infer that it contains $Z(G)$. The second statement follows from Corollary 4.5.5 of \cite{18}.
\item[{\rm c)}] Suppose $H\leq G$ with $|H|=p$ and $H\neq Z(G)$. Then there is $K\in{\cal CD}(G)$ such that $|K|=p^2$ and $H\subset K$. We get $C_G(K)\subseteq C_G(H)$ and since $C_G(K)$ is of order $p^{n-1}$, we obtain $C_G(H)=C_G(K)$. Thus $m_G(H)=|H||C_G(H)|=p\cdot p^{n-1}=p^n$. 
    
    Now, suppose $A\leq G$ with $|A|\geq p^2$ and $A\notin{\cal CD}(G)$. Then there are $B,C\in{\cal CD}(G)$ such that $B\subset A\subset C$ and $|A|/|B|=|C|/|A|=p$. It follows that $C_G(C)\subseteq C_G(A)\subseteq C_G(B)$ and $|C_G(B)|/|C_G(C)|=p^2$. We cannot have $C_G(A)=C_G(B)$ because this implies $m_G(A)>m^*(G)$, a contradiction. Also, we cannot have $|C_G(A)|=p|C_G(C)|$ because this implies $A\in{\cal CD}(G)$, a contradiction. So, $C_G(A)=C_G(C)$, which leads to $m_G(A)=p^n$. Finally, note that $m_G(1) = |G| = p^n$.  Thus ${\rm Im}(m_G)=\{p^n,p^{n+1}\}$.\qed
\end{itemize}

\bigskip\noindent{\bf Proof of Theorem 1.2.} 
If $G$ is a nonabelian group of order $pq$, $p$ and $q$ primes, then $1 < N < G$ where $N$ is a normal subgroup of prime order, and ${\cal CD}(G) = \{ N \}$, and thus $G$ has dense ${\cal CD}$-subgroups.

Conversely, suppose that $G$ has dense ${\cal CD}$-subgroups.
Let $|G|=p_1^{n_1}\cdots p_k^{n_k}$ with $k\geq 2$. If $n_i=1$, $\forall\, i=1,...,k$, then $G$ is a ZM-group and so its Chermak--Delgado lattice is a chain of length $0$ by Theorem 3 of \cite{12}. This satisfies the density property if and only if $k=2$, that is, $G$ is a nonabelian group of order $p_1p_2$.

In what follows, we assume that $n_1\geq 2$ and let $H_1\leq G$ with $|H_1|=p_1^2$. Then between $1$ and $H_1$ there is a subgroup $K_1\in{\cal CD}(G)$ of order $p_1$. Since $Z(G)\subseteq K_1$, we have the following two cases:

\bigskip\hspace{5mm}{\bf Case 1.} $Z(G)=K_1$.

\noindent Then $n_2=\dots=n_k=1$ and so $|G|=p_1^{n_1}p_2\cdots p_k$. Indeed, if $n_i\geq 2$ for some $i=2,...,k$, then there is $K_i\in{\cal CD}(G)$ of order $p_i$. This implies $Z(G)\subseteq K_i$, that is $K_1\subseteq K_i$, which leads to $p_1\mid p_i$, a contradiction.

Clearly, we can choose $M_1\leq G$ with $|M_1|=p_1^{n_1-1}$ and $M_1\in{\cal CD}(G)$.  Now, subgroups in ${\cal CD}(G)$ are permutable, and ${\cal CD}(G)$ is closed under conjugation, and so the normal closure $M_1^G$ is a $p_1$-subgroup.  Now the quotient group $G/M_1^G$ is a ZM-group and hence solvable, and $M_1^G$ is a $p_1$-subgroup which is solvable, and so it follows that $G$ is solvable.  (Note: the solvability of a ZM-group can be proven as a consequence of Burnside's normal complement theorem.  Even more is true: in fact a ZM-group is metacyclic.)

We observe that $k=2$. Indeed, if $k\geq 3$, then as $G$ is solvable it contains a Hall subgroup $H_{23}$ of order $p_2p_3$. Let $K\leq G$ such that $1\subset K\subset H_{23}$ and $K\in{\cal CD}(G)$. Then $|K|\in\{p_2,p_3\}$ and $Z(G)\subseteq K$, implying that either $p_1\mid p_2$ or $p_1\mid p_3$, a contradiction. Thus
\begin{equation}
|G|=p_1^{n_1}p_2 \mbox{ and } m^*(G)=p_1^{n_1+1}p_2.\nonumber
\end{equation}

So $|C_G(M_1)|=p_1^2p_2$ and let $P_2$ be a Sylow $p_2$-subgroup of $G$ that centralizes $M_1$.  Now, the interval $(M_1,G)$ of $L(G)$ contains ${\cal CD}$-subgroups and such a subgroup is either a Sylow $p_1$-subgroup $P_1$ of $G$ or a product $M_1 P_2^y$ for some $y \in G$.  If $P_1\in{\cal CD}(G)$, we get $|C_G(P_1)|=p_1p_2$ and therefore there is $x\in G$ such that $P_2^x$ centralizes $P_1$. Since $P_2^x$ centralizes $P_2^x$, it follows that $P_2^x$ centralizes $P_1P_2^x=G$. So, $P_2^x\subseteq Z(G)$, implying that $p_2\mid p_1$, a contradiction. We now argue that $M_1$ is normal in $G$.  Otherwise, the normal closure $M_1^G$ is a Sylow $p_1$-subgroup of $G$ that is in ${\cal CD}(G)$, a contradiction. So $M_1$ is normal in $G$ and suppose $M_1P_2^y\in{\cal CD}(G)$ for some $y \in G$.  Then $M_1P_2\in{\cal CD}(G)$ and so $|C_G(M_1P_2)|=p_1^2$. But $P_2$ centralizes $M_1P_2$ and consequently $p_2\mid p_1^2$, a contradiction.

Thus this case is impossible.

\bigskip\hspace{5mm}{\bf Case 2.} $Z(G)=1$.

\noindent Since $M(G)\subseteq K_1$, we have either $M(G)=1$ or $M(G)=K_1$. 

If $M(G)=1$ then $1\in{\cal CD}(G)$ and so ${\cal CD}(G)$ does not contain any subgroup of prime order by Corollary 7 of \cite{11}. This contradicts our hypothesis.

Assume now that $M(G)=K_1$. Similarly with Case 1, we get $n_2=\dots=n_k=1$ and $k=2$, that is $|G|=p_1^{n_1}p_2$. Since $|C_G(K_1)|$ is a proper divisor of $|G|$ and $p_1|C_G(K_1)|=m_G(K_1)>m_G(1)=|G|$, it follows that $|C_G(K_1)|=p_1^{n_1}$ and hence $m^*(G)=p_1^{n_1+1}$.  Let $P_2$ be a Sylow $p_2$-subgroup of $G$.  Since $|G|=p_1^{n_1}p_2$ and $n_1 \geq 2$, we have $P_2 < K_1P_2 < G$, and so $P_2$ is not maximal in $G$.  Thus, there is $X \in{\cal CD}(G)$ so that $P_2 < X < G$, and hence $p_2 \mid |X|$, contradicting $m^*(G)=p_1^{n_1+1}$.  \qed

\vspace*{3ex}
\small

\begin{minipage}[t]{8cm}
Ryan McCulloch\\
Department of Mathematics and Statistics\\
Binghamton University\\
Binghamton, NY 13902\\
e-mail: {\tt rmccullo1985@gmail.com}
\end{minipage}
\hfill
\begin{minipage}[t]{5cm}
Marius T\u arn\u auceanu \\
Faculty of  Mathematics \\
``Al.I. Cuza'' University \\
Ia\c si, Romania \\
e-mail: {\tt tarnauc@uaic.ro}
\end{minipage}

\end{document}